\documentclass{article}

\usepackage{graphicx,color,amsmath,amssymb}
\usepackage{fancybox,%
            epsf,%
            theorem,%
            eepic,%
            epsfig,%
            amsfonts,%
            pifont}%
\usepackage{graphics}
\usepackage{color}

\newcommand{\Tau}{{\cal T}}

\newcommand{\Id}{\mathrm{Id}}

\newcommand{\mR}{\mathbb R}
\newcommand{\mRp}{\mathbb R_{\geq 0}}
\newcommand{\cKL}{\mathcal{KL}}

\newcommand{\norm}[1]{\left|{#1}\right|}

\newcommand{\la}{\left\langle}
\newcommand{\ra}{\right\rangle}

\renewcommand{\r}[1]{(\ref{#1})}
\newcommand{\Ar}{\mathbb{R}}
\newcommand{\N}{\mathbb{N}}

\newcommand{{\R}}{\Ar}
\newcommand{\beq}{\begin{equation}}
\newcommand{\eeq}{\end{equation}}
\newcommand{\be}{\begin{equation}}
\newcommand{\ee}{\end{equation}}
\newcommand{\bea}{\begin{eqnarray}}
\newcommand{\eea}{\end{eqnarray}}
\newcommand{\br}{\begin{eqnarray}}
\newcommand{\er}{\end{eqnarray}}
\newcommand{\brs}{\begin{eqnarray*}}
\newcommand{\ers}{\end{eqnarray*}}
\newcommand{\ba}{\begin{array}}
\newcommand{\ea}{\end{array}}

\newcommand{\caG}{{\cal {G}}}

\newcommand{\lb}{\lambda}

\newcommand{\al}{\alpha}

\newcommand{\lp}{\left(}
\newcommand{\rp}{\right)}

\newcommand{\bed}{\begin{description}}
\newcommand{\eed}{\end{description}}

\def\cf{{\it cf.}}
\def\n{\noindent}
\def\rref#1{(\ref{#1})}
\def\eps{\varepsilon}

\def\al{\alpha}

\def\th{\theta}

\def\mt{\mapsto}
\def\EOP{\ \hfill \rule{0.5em}{0.5em} }

\newtheorem{thm}{Theorem}
\newtheorem{lemma}[thm]{Lemma}

\newtheorem{prop}[thm]{Proposition}
\begin{theorembodyfont}{\rmfamily}
\newtheorem{rmk}[thm]{Remark}

\end{theorembodyfont}

\newtheorem{deff}[thm]{Definition}

\begin{document}

\title{\bf \Large \vspace{5mm}%
Uniform stabilization for linear systems with persistency of excitation. The
neutrally stable and the double integrator cases.}
\author{Antoine Chaillet, \qquad
        Yacine Chitour, \qquad
        Antonio Lor\'{\i}a,\qquad
        Mario Sigalotti}
\date{}
\sloppy
\maketitle
\parskip=1pt plus 1pt minus 1pt
\global\setlength\theorempreskipamount{2pt} 
\global\setlength\theorempostskipamount{2pt} 

\begin{abstract}
\let\piedepag=\thefootnote\def\thefootnote{}
\phantom{a}\footnote{{\indent  A. Chaillet is with Centro di Ricerca Piaggio,
Facolt\`a di ingegneria, Via Diotisalvi 2, 56100 Pisa, Italy; Y. Chitour and
 A. Lor\'{\i}a are with Laboratoire des Signaux et Syst\`emes, Sup\'elec,
 3, Rue Joliot Curie,  91192 Gif s/Yvette, France; Y. Chitour is also with
Universit\'e Paris Sud, Orsay and A. Lor\'{\i}a is with C.N.R.S;
M. Sigalotti is with INRIA, Institut \'Elie
Cartan, UMR 7502 Nancy-Universit\'e/CNRS/INRIA,
 POB 239, Vand\oe uvre-l\`es-Nancy 54506,  France.
E-mails:  {\tt chaillet@ieee.org},
{\tt chitour@lss.supelec.fr}, {\tt loria@lss.supelec.fr}, {\tt mario.sigalotti@inria.fr}.
}}
\let\thefootnote=\piedepag
\mbox{}\hspace{-3mm}%
Consider the controlled system $dx/dt = Ax + \alpha(t)Bu$ where the
pair $(A,B)$ is stabilizable
and $\alpha(t)$ takes values in $[0,1]$
and is persistently exciting, i.e., there exist two positive constants $\mu,T$
such that, for every $t\geq 0$, $\int_t^{t+T}\alpha(s)ds\geq \mu$.
In particular, when $\alpha(t)$
becomes zero the system dynamics switches to an uncontrollable
system. In this paper, we address the following question:
is it possible to
find a linear time-invariant state-feedback $u=Kx$, with
$K$ only depending on $(A,B)$ and possibly on $\mu,T$,
which globally asymptotically stabilizes the system? We
give a positive answer to this question for two cases: when $A$ is
neutrally stable and when the system is the double integrator.\\
\end{abstract}

\noindent {\bf Notation.}
 A continuous function $\phi : \mRp \to \mRp$ is of \emph{class
${\cal K}$} ($\phi \in {\cal K}$), if it is strictly increasing
\textcolor{black}{and $\phi(0)=0$}.
$\psi : \mRp\to\mRp$ is of class $\cal L$ ($\psi\in\cal L$) if it
is continuous, non-increasing and tends to zero as its argument tends to
infinity. A function $\beta: \mRp\times \mRp \to \mRp$ is said to be
a \emph{class $\cal{KL}$-function} if, $\beta(\cdot,t)\in\cal K$ for
any $t\geq 0$, and $\beta(s,\cdot)\in\cal L$ for any $s\geq 0$. We
use $\norm{\cdot}$ for the Euclidean norm of vectors and the induced
$L_2$-norm of matrices.

\section{Introduction}

In the paper \cite{DOUBLEINT} we posed the following problem: consider the system
$\dot x = f(t,x) + u$ with $x\in\mR^n$ and the stabilizing control $u=u^*$. Consider
now the system
\begin{equation}
\label{1000} \dot x = f(t,x) + \alpha(t)u \,,
\end{equation}
where rank$\{\alpha(t)\} \neq n$ for certain times $t$
({\it i.e.} $\alpha$ may be rank-deficient over possibly ``large'' intervals of time).
Under which conditions
imposed on $\alpha$ does the closed-loop system  \rref{1000} with
the {\em same} control $u^*$ is asymptotically stable?
It must be stressed that a complete knowledge of
$\alpha$ (and, in particular, precise information on the set of times
where it is rank deficient) would be a too restrictive condition to
impose on $\alpha$. We are rather looking for a condition valid for a
whole class $\cal {G}$ of functions $\alpha$ and, therefore, we
expect the closed-loop systems \rref{1000} with $u^*$ to be
asymptotically stable for every $\alpha\in {\cal {G}}$.

In order to characterize such a condition, let us consider a similar
problem stemming from identification and adaptive control. It concerns
the linear system $\dot x = -P(t)u$, where the matrix
$P(\cdot)$ is symmetric non-negative and now plays the role of $\alpha$. If $P\equiv I$
then $u^*=x$ stabilizes the system exponentially. But what if $P(t)$ is only
semidefinite for all $t$? \textcolor{black}{Under which conditions}
does $u^*=x$ still stabilize the system?  For this
particular case the answer to this question can be found in the literature:
from the seminal paper \cite{MORNAR} we know that for the system
\begin{equation}
\label{1001}
\dot x = -P(t)x
\end{equation}
with $x\in\mR^n$, $P\geq0$ bounded and with bounded derivative, it is
{\em necessary and sufficient}, for global exponential stability, that $P$ be also
{\em persistently exciting} (PE), {\em i.e.}, that there exist $\mu>0$ and $T>0$ such that
\begin{equation}\label{1002}
\int_{t}^{t+T} \xi^\top P(\tau)\xi \geq \mu
\end{equation}
for all unitary vectors $\xi\in\mR^n$ and all $t\geq 0$. Therefore, as regards the
stabilization of \rref{1000}, the notion of
{\em persistent excitation} seems to be a reasonable additional assumption to consider
for the signals $\alpha$.

In this paper, we focus on $n$-dimensional linear time-invariant systems
\begin{equation}\label{system}
\dot x = A x +\alpha(t)Bu\,,
\end{equation}
where the input perturbation $\alpha$ is a {\em scalar} PE signal,
i.e., $\alpha$ takes values in $[0,1]$ and
there exist
two positive constants $\mu,T$ such that, for every $t\geq 0$,
\begin{equation}
\label{Tmu}
\int_t^{t+T}\alpha(s)ds\geq \mu.
\end{equation}
Given two positive real numbers $\mu\leq T$, we denote by ${\cal
{G}}(T,\mu)$ the class of all PE signals verifying \rref{Tmu}. Note
that we do not consider here any extra assumption on the regularity
of the PE signal $\alpha$ (e.g. having a bounded derivative).


An interpretation of the stabilization mechanism can be given, in the case of scalar
systems of type \rref{1001}, in terms of ``average''. Roughly speaking, one can dare
say that, even though
it is not the control action $u^*$ that enters the system for each $t$, this
``ideal'' control does drive the system ``{\em in average}''.
Indeed, for $t\geq 0$, set
$\la\alpha\ra\!\!(t):=\frac1T\int_t^{t+T}\alpha(s)ds$ for any PE signal $\alpha$.
Then, along non-trivial trajectories of $\dot x(t)=-\alpha(t) x(t)$, one has
for every $t\geq 0$,
$$
\frac1T\ln\lp\frac{x(t+T)}{x(t)}\rp=\frac1T\int_t^{t+T}\frac{\dot x(s)}{x(s)}ds
=-\la\alpha\ra(t).
$$
The control action $u=-\alpha(t)x$ which, in
average, corresponds
to $u=-\la\alpha\ra\!\!(t)x$, is tantamount to applying $u^*=-x$, modulo a gain-scale that
only affects the rate of convergence but not the stabilization property of $u^*$.
Of course the previous naive thinking largely relies on the fact that we are dealing
with a scalar system. Indeed, persistency of excitation does not guarantee the existence
of an averaged system in the sense of {\it e.g.} \cite{TEPEAE99}. This observation makes
all the related techniques unapplicable.

Our main goal in this paper consists of stabilizing \rref{system} to the origin
with a linear feedback. Therefore, we will assume in the sequel that
$(A,B)$ is a stabilizable pair. It is obvious that if $\al\equiv 1$, for
a proper choice of $K$, we have that the control $u=u^*$ with $u^*=
Kx$ renders the closed-loop system globally exponentially stable. If
$\alpha$ is not constant, consider the following question:
\begin{quote}
\label{ques:01} (Q1-0)\ \ Does  $u=\alpha(t)u^*$, with $\alpha$ an
arbitrary PE signal, stabilize \rref{system}?
\end{quote}
For systems of the type \rref{system}, the answer to Question
($Q$1-0) is negative in general. Indeed, the scalar case
essentially corresponds to stabilizing
$\dot x(t)= \lambda x(t)+\alpha(t)u$. Using the linear feedback
$u^*=k x$ for some $k<0$, one gets, after simple computations, that
for $t\geq 0$,
\begin{equation}\label{triv}
\frac1T\ln\lp\frac{x(t+T)}{x(t)}\rp=\lambda+\la\alpha\ra\!(t)k\leq
\lambda+\frac{\mu}Tk,
\end{equation}
if $\alpha$ is a PE-signal verifying \rref{Tmu} for some fixed positive
constants $\mu\leq T$. One deduces from
\rref{triv} that global exponential stabilization occurs if the
negative constant $k$ is chosen so that $k<-\frac T{\mu}\lambda$. If
$\lambda>0$, then the choice of $u^*=k x$ depends on $\mu,T$, the parameters
of the PE-signal. Therefore, question ($Q$1-0) \textcolor{black}{may only receive a
positive answer} for systems \rref{system} with $A$ marginally stable, {\it i.e.} all the
eigenvalues of $A$ have non-positive real part:
\begin{quote}
\label{ques:1} (Q1)\ \ Given $A$ marginally stable, does  $u=\alpha(t)u^*$, with $\alpha$
an arbitrary PE signal, stabilize \rref{system}?
\end{quote}
Intuitively one may think that the global exponential stability of
the closed-loop system is guaranteed, at least, for a particular choice of
$K$ and a particular class of PE functions
$\alpha$. In that spirit, consider the following question:
\begin{quote}
\label{ques:2} (Q2)\ \ Given a class ${\cal {G}}(T,\mu)$ of PE signals, can one
determine a linear feedback $u^*= Kx$ such
that $u=\al(t)u^*$ stabilizes \rref{system}, \textcolor{black}{for all PE signal
$\alpha$ in} ${\cal {G}}(T,\mu)$?
\end{quote}
While, for the scalar case, the answer is clearly positive (as shown previously), the
general case $n>1$ is fundamentally different and, in view of the
available tools from the literature of adaptive control, a proof (or
disproof) of the conjecture above is far from evident. A first step
in the solution of $(Q2)$ for the case $n=2$ has been undertaken
in \cite{DOUBLEINT} where we showed that, for certain persistently
exciting functions $\alpha$ and certain values of $K
$,
exponential stability follows. \textcolor{black}{We underline the fact that,
in Question $(Q2)$, the gain $K$ is required to be valid \emph{for all} signals
in the considered class.}

A particular case of this problem may be re-casted in the context of
switched systems --\cf\, \cite{LIBBOOK}. Indeed, consider the particular case
in which the PE signal $\alpha$ takes only the values $0$ and $1$.
In this setting, the system \rref{system} after the choice
$u=u^*=Kx$, switches between the uncontrolled system $\dot x = Ax$
and the exponentially stable system
$\dot x= (A+BK)x$. At this point,
it is worth emphasizing that since switches occur between a possibly
unstable dynamics and a stable one, Lyapunov-based conditions for
switched system between stable dynamics are inapplicable; see for
instance \cite{SHOWIRMASWULKIN} for results using common quadratic
Lyapunov functions, \cite{BRANICKY98,COLGERAST} for theorems
relying on multiple Lyapunov functions and also \cite{BOSCAIN2002} for a more geometric
approach.


The first issue we address in this paper regards the
controllability of \rref{system} uniformly with respect to $\alpha\in {\cal {G}}(T,\mu)$.
If the pair $(A,B)$ is controllable, we prove that \rref{system} is
(completely) controllable in time $t$ if and only if $t>T-\mu$.

We next focus on the stabilization of \rref{system} by a linear
feedback $u=Kx$, that is, we address Question $(Q2)$ for system
\rref{system}. Namely, we look for the existence of a matrix $K$ of
size $m\times n$ such that, for every $\alpha\in {\cal {G}}(T,\mu)$, the origin of the
system
\[
\dot x=(A+\alpha(t) BK) x
\]
is globally asymptotically stable. It is
of course assumed that $(A,B)$ is stabilizable. We first treat the case where $A$ is
neutrally stable. We actually
determine a feedback $K$ as required, which in addition provides a
positive answer to Question $(Q1)$ posed for system \rref{system}.

Finally, we consider the case of the double integrator , \textcolor{black}{that is,
$$
A=\begin{pmatrix}
0&1\\
0&0
\end{pmatrix}, \ \ \
B=\begin{pmatrix} 0\\1
\end{pmatrix}\,.
$$
}

In \cite{DOUBLEINT} we
already studied such a system  under the assumption that
$\alpha$ is PE. The solutions given in that paper, however, do not
bring a complete answer to the questions posed above:
the first solution relies on backstepping and, therefore,
requires a bound on the derivative of $\alpha$ while the second is based on
a normalization of $\alpha$ \textcolor{black}{and imposes a relationship between the
parameters $T$ and $\mu$ involved in the persistency of excitation.} In the present
paper, we bring a positive
answer to Question $(Q2)$ for every class ${\cal {G}}(T,\mu)$ of PE signals
and a negative answer to Question \textcolor{black}{($Q$1)}.

%

The rest of the paper is organized as follows. In coming section we
provide the main notations and the result on the controllability of
multi-input linear systems subject to PE signals. We discuss
stabilization issues
in Section~\ref{neutr} for $A$
neutrally stable and in Section~\ref{sec:mainproof} for the double
integrator, \textcolor{black}{first in case of a scalar control
and then in a more general setup}. We
close the paper with an Appendix which
contains the proof of a crucial
technical
result.

\section{Notations and basic result on controllability}
In this paper, we are concerned with linear systems subject to a
scalar persistently exciting signal, {\em i.e.}, systems of the type \rref{system}
where $\alpha$ is a PE signal. The latter is defined as follows.

\begin{deff}[$(T,\mu)$-signal]\label{Tm-signal}
Let $\mu\leq T$ be positive constants. A \emph{$(T,\mu)$-signal} is
a measurable function $\alpha:{\mathbb R_{\geq 0}}\to[0,1]$
satisfying
\begin{equation}\label{EP}
\int_t^{t+T}\alpha(s)ds \geq \mu\,,\quad \forall t\in\mRp\,.
\end{equation}
We use ${\cal {G}}(T,\mu)$ to denote the set of all $(T,\mu)$-signals.
\end{deff}

\textcolor{black}{Notice that, for any such signal $\alpha$, existence and uniqueness
of the solutions of \rref{system} is guaranteed.}

\begin{rmk}\label{notime}
If $\alpha(\cdot)$ is a $(T,\mu)$-signal, then, for every $t_0\geq
0$, $\alpha(t_0+\cdot)$ is again a $(T,\mu)$-signal.
\end{rmk}

\textcolor{black}{Our first result studies the following property for \rref{system}.}

\begin{deff}[Controllability in time $t$]
We say that system \rref{system} is controllable in time $t>0$ for
${\cal {G}}(T,\mu)$ if, for every $\alpha\in{\cal {G}}(T,\mu)$, the
time-varying linear controlled system defined by \rref{system} is
completely controllable in time $t$.
\end{deff}

More precisely, we establish the following result on the
controllability of system \rref{system}.

\begin{prop}\label{cont-t}
Let $\mu\leq T$ be two positive constants and $(A,B)$ a controllable
pair of matrices of size $n\times n$ and $n\times m$ respectively.
Then, system \rref{system} is controllable
in time $t$ for ${\cal {G}}(T,\mu)$ if and only if $t>T-\mu$.
\end{prop}

\underline{Proof of Proposition~\ref{cont-t}.} Following the
classical proof of the Kalman condition for the controllability of
autonomous linear systems \textcolor{black}{(\cf {\it e.g} \cite{SONTAGBOOK})},
it is easy to see that the conclusion
does not hold if and only if there exists a non-zero vector $p\in
\R^n$ and a $(T,\mu)$-signal $\alpha$ such that the function
\begin{equation}\label{Kal}
\alpha(s)p^\top e^{A(t-s)}B= 0, \ \ \hbox{a.e. in } [0,t].
\end{equation}
If $t>T-\mu$, there exists a subset $J$ of $[0,T]$
of positive measure such that $\alpha(s)>0$ for $s\in J$. Therefore,
the real-analytic function $s\mt p^\top e^{A(t-s)}B$ is equal to
zero on $J\subset [0,T]$. It must then be identically equal to zero,
which implies that $(A,B)$ is not controllable. We reach a
contradiction and one part of the equivalence is proved.
If $t\leq T-\mu$, ${\cal {G}}(T,\mu)$
contains a PE-signal identically equal to zero on time intervals of
lengtht $t$ and any non-zero vector $p$ verifies \rref{Kal}.\EOP

\textcolor{black}{The rest of the paper is concerned with the stabilization of
\rref{system}.
We address} the following problem. Given $T\geq \mu>0$, we want to find a matrix $K$
of size $m\times n$ which makes the origin of
\begin{equation}\label{feedback}
\dot x=(A+\alpha(t) BK) x
\end{equation}
globally asymptotically stable, uniformly with respect to every
$(T,\mu)$-signal $\alpha$ \textcolor{black}{({\it i.e.}, $K$ is required to depend only
on $A$, $B$, $T$ and $\mu$ and to be valid for all signals $\alpha$ in the class
$\mathcal G(T,\mu)$)}. Referring to
$$x(\cdot\,;t_0,x_0,K,\alpha)=
(x_1(\cdot\,;t_0,x_0,K,\alpha),\dots,x_n(\cdot\,;t_0,x_0,K,\alpha))^\top,
$$
as the solution of \rref{feedback} with initial condition
$x(t_0;t_0,x_0,K,\alpha)=x_0$, we introduce the following
definition.

\begin{deff}[$(T,\mu)$-stabilizer]\label{stab}
Let $\mu\leq T$ be positive constants. The gain $K$ is said
to be a \emph{$(T,\mu)$-stabilizer} for \rref{system} if there
exists a class $\cKL$-function $\beta$ such that, for every
$\al\in\caG(T,\mu)$, $x_0\in\mR^n$, and $t_0\in\mRp$, the
solution of \rref{feedback} satisfies
$$
\norm{x(t;t_0,x_0,K,\alpha)}\leq \beta(\norm{x_0},t-t_0)\,,\quad
\forall t\geq t_0\,.
$$
\end{deff}
\begin{rmk}\label{exp0}
Since we are dealing with linear systems (in the state), it is a
standard fact that one can rephrase the above definition as follows:
the gain $K$ is said to be a \emph{$(T,\mu)$-stabilizer} for
\rref{system} if \rref{feedback} is exponentially stable, uniformly
with respect to every $(T,\mu)$-signal $\alpha$. Here, uniformity
means that the rate of (exponential) decrease only depends on
$(A,B)$ and $\mu,T$.
\end{rmk}

\section{The neutrally stable case}\label{neutr}
The purpose of this section consists of proving the following
theorem.
\begin{thm}\label{thm-neutrally}
Assume that the pair $(A,B)$ is stabilizable and that the matrix $A$
is neutrally stable.
Then there exists a matrix $K$ of size $m\times n$ such that, for
every $0<\mu\leq T$, the gain $K$ is a \emph{$(T,\mu)$-stabilizer}
for \rref{system}.
\end{thm}

Since the feedback $K$ determined above does not depend on the
particular class ${\cal {G}}(T,\mu)$, we bring a positive answer to
Question $(Q1)$ in the case where $A$ is neutrally stable.
\begin{rmk}
It can be seen along the proof below that
the gain $K=-rB^\top$, with an arbitrary $r>0$, does the job in
the case where $A$ is skew-symmetric and $(A,B)$ is controllable.
\end{rmk}



\underline{Proof of Theorem~\ref{thm-neutrally}.}  Recall that a
matrix is said to be neutrally stable if its eigenvalues have
non-positive real part and those with zero real part have trivial
corresponding Jordan blocks. The proof of
Theorem~\ref{thm-neutrally} is based on the following equivalence
result.

\begin{lemma}\label{reduction}
It is enough to prove Theorem~\ref{thm-neutrally} in the case where $A$ is skew-symmetric
and $(A,B)$ controllable.
\end{lemma}

\underline{Proof of Lemma~\ref{reduction}.} Let $(A,B)$ be a
stabilizable pair with $A$ neutrally stable. Since the non-controlled
part of the linear system $\dot x=Ax+Bu$ is already stable, it is
enough to focus on the controllable part of $(A,B)$. Hence, we
assume that $(A,B)$ is controllable. Up to a linear change of
variable, $A$ and $B$ can be written as
$$
A=\begin{pmatrix}
A_1&A_2\\
0&A_3
\end{pmatrix}, \ \ \
B=\begin{pmatrix} B_1\\B_3
\end{pmatrix},
$$
where $A_1$ is Hurwitz and all the eigenvalues of $A_3$ have zero
real part. From the neutral stability assumption, we deduce that
$A_3$ is similar to a skew-symmetric matrix. Up to a further linear change of
coordinates, we may assume that $A_3$ is indeed skew-symmetric. From the
controllability assumption, we deduce that $(A_3,B_3)$ is
controllable. Setting $x=(x_1^\top,x_3^\top)^\top$ according to the above decomposition,
the system \rref{system} can be written as
\br
\dot x_1&=&A_1x_1+A_2x_3+\alpha(t)B_1u,\label{h1}\\
\dot x_3&=&A_3x_3+\alpha(t)B_3u.\label{h2}
\er
Assume that Theorem~\ref{thm-neutrally} holds for \r{h2}, {\em i.e.}, there
exists $K_3$ such that, for
every $0<\mu\leq T$, the gain $K_3$ is a \emph{$(T,\mu)$-stabilizer}
for \r{h2}. Take
$$K=\lp\ba{c}{\bf 0}\\ K_3\ea\rp,$$
where all entries of ${\bf 0}$ are null. Then the conclusion for
system \rref{system} follows, since an autonomous linear Hurwitz
system subject to a perturbation whose norm converges exponentially
(with respect to time) to zero is still asymptotically stable at the
origin. \EOP


\vspace{0.5cm} Based on Lemma \ref{reduction}, we assume that $A$ is
skew-symmetric and $(A,B)$ is controllable for the rest of the argument.
We will prove that the
gain $K=-B^\top$ does the job.

We consider the time derivative of the Lyapunov function $V(x)=\norm x^2/2$
along non-trivial solutions of the closed-loop system
\begin{equation}\label{feedback2}
\dot x=\left(A-\alpha(t)BB^\top \right)x\,,
\end{equation}
and get $\dot V=-\alpha(t)\norm{B^\top x}^2$, which can be also written
$$
\frac{\dot V}{V}=-\alpha(t)\frac{
\norm{B^\top x}^2}{V}\,.
$$
Integrating both parts and defining $x(\cdot)$ as
$x(\cdot\,;t_0,x_0,K,\alpha)$ and $v(\cdot)=V(x(\cdot))$, where
$(t_0,x_0)$ and $\al$ denote, respectively, an arbitrary initial condition and an
arbitrary $(T,\mu)$-signal given arbitrary $T\geq\mu>0$, we get that
\begin{equation}\label{Vint}
\int_{t_0}^{t_0+T} \frac{\dot v(t)}{v(t)}dt=-\int_{t_0}^{t_0+T}
\alpha(t)\frac{\norm{B^\top x(t)}^2
}{v(t)}dt\,.
\end{equation}

\begin{lemma}\label{claim1}
For every $0<\mu\leq T$, there exists a positive constant $\eta$
such that, for any $(T,\mu)$-signal $\alpha$ and any initial state
$\norm{x_0}=1$, it holds that
$$
\int_{t_0}^{T+t_0} \alpha(t)\frac{\norm{B^\top x(t)}^2
}{v(t)}dt \geq
\eta\,.
$$
\end{lemma}

\underline{Proof of Lemma~\ref{claim1}.}  Because of
Remark~\ref{notime} we take, without loss of generality, $t_0=0$. We fix
$0<\mu\leq T$ and reason by contradiction, {\it i.e.}, we assume
that there exist a sequence $\{x_{0i}\}_{i\in\N}$ such that
$\norm{x_{0i}}=1$ for all $i\in\N$ and a sequence of
$(T,\mu)$-signals $\alpha_i$ such that
\begin{equation}\label{limint}
\lim_{i\to\infty} \int_{0}^{T} \alpha_i(t)\frac{\norm{B^\top
x_i(t)}^2}{v_i(t)}dt =0\,,
\end{equation}
where $x_i(\cdot)$ denotes $x(\cdot\,;0,x_{0i},K,\alpha_i)$ and
$v_i(\cdot)=V(x_i(\cdot))$. Since $\{x_{0i}\}_{i\in\N}$ belongs to a compact set,
there exists a
subsequence $\{x_{0i_j}\}_{j\in\N}$ such that
$$
\lim_{j\to\infty} x_{0i_j}=x_{0\star}\,,\quad \textrm{with}\quad
\norm{x_{0\star}}=1\,.
$$
On the other hand,
recall that the space $L^\infty(\mRp,[0,1])$ is sequentially weakly-$\star$
compact (see, for instance, \cite[Chapter IV]{brezis83}), that is,
for every sequence $\{\beta_i\}_{i\in\mathbb N}\subset L^\infty(\mRp,[0,1])$
there exist $\beta_\star\in L^\infty(\mRp,[0,1])$ and a subsequence
$\{\beta_{i_j}\}_{j\in\mathbb N}$ such that
for every $\varphi\in L^1(\mRp,\mR)$ the following holds
\begin{equation}\label{weakstar}
\lim_{j\to\infty} \int_0^\infty \beta_{i_j}(s)\varphi(s)ds
=      \int_0^\infty \beta_\star(s) \varphi(s)ds\,.
\end{equation}

Therefore, we can extract a
subsequence $\{\alpha_{i'_j}\}_{j\in{\mathbb N}}$ of
$\{\alpha_{i_j}\}_{j\in{\mathbb N}}$ which
converges weakly-$\star$ to a measurable
function $\alpha_\star$. Note that $\alpha_\star$
is itself a
$(T,\mu)$-signal.

The convergence of $\alpha_{i'_j}$ to $\alpha_\star$ implies a convergence
of the solutions
of their corresponding non-autonomous linear dynamical systems in the form
\r{feedback2}: this is a particular case of the classical Gihman
convergence with respect to
time-varying parameters for ordinary differential equations (see \cite{gihman}
and also \cite{liu-sus} for a general discussion in the framework of control
theory). For the sake of completeness, and because the situation faced here
can be handled  with a specific and simpler approach, we include in the
Appendix the proof of Proposition~\ref{technic}, whose statement allows to conclude that
$\{x_{i'_j}(\cdot)\}_{j\in\mathbb N}$ converges, uniformly on compact
time intervals, to
$x_{\star}(\cdot):=x(\cdot\,;0,x_{0\star},K,\alpha_\star)$ as $j$
tends to infinity. Letting $v_\star(\cdot)=V(x_\star(\cdot))$, one has,
for every $j\in \mathbb N$,
$$
\int_{0}^{T} \alpha_{i'_j}(t)\frac{\norm{B^\top x_{i'_j}(t)^2}}{v_{i'_j}(t)}dt
-\int_{0}^{T} \alpha_\star(t)\frac{\norm{B^\top
x_\star(t)}^2}{v_\star(t)}dt=
$$
$$
\int_{0}^{T} \big(\alpha_{i'_j}(t)-
\alpha_\star(t)\big)\frac{\norm{B^\top
x_\star(t)}^2}{v_\star(t)}dt+
\int_{0}^{T} \alpha_{i_j'}(t)\lp\frac{\norm{B^\top x_{i'_j}(t)^2}}{v_{i'_j}(t)}-
\frac{\norm{B^\top
x_\star(t)}^2}{v_\star(t)}\rp dt.
$$
Letting $j$ tend to infinity, taking into account the weak-$\star$ convergence of
$\alpha_{i'_j}$ to $\alpha_\star$,
 and applying \r{limint},
we get that
\brs
\int_{0}^{T} \alpha_\star(t)\frac{\norm{B^\top
x_\star(t)}^2}{v_\star(t)}dt &\leq& \liminf_{j\to\infty}\int_{0}^{T} \alpha_{i'_j}(t)\left|\frac{\norm{B^\top x_{i'_j}(t)^2}}{v_{i'_j}(t)}-
\frac{\norm{B^\top
x_\star(t)}^2}{v_\star(t)}\right|dt\\
&\leq& \liminf_{j\to\infty}\int_{0}^{T} \left|\frac{\norm{B^\top x_{i'_j}(t)^2}}{v_{i'_j}(t)}-
\frac{\norm{B^\top
x_\star(t)}^2}{v_\star(t)}\right|dt\\
&=&0\,,
\ers
where the last equality follows from the uniform convergence of
$x_{i'_j}(\cdot)$ to $x_\star(\cdot)$ on the interval $[0,T]$.
Hence, we conclude that
\begin{equation}\label{zeroae}
\alpha_\star(t)B^\top x_\star(t)=0
\end{equation}
for almost every $t\in[0,T]$.
Furthermore, since $x_{\star}$ is a solution of \rref{feedback2}, we
get from \rref{zeroae} that $x_\star(t)=e^{At}x_{0\star}$ for
$t\in[0,T]$. Moreover, since $\alpha_{\star}$ is a $(T,\mu)$-signal,
it is strictly positive on a subset $J$ of $[0,T]$ of positive
measure, and then, according to \rref{zeroae}, $B^\top e^{At}x_{0\star}$ must be
equal to zero on $J$. Since the latter is real-analytic, it must be
identically equal to zero on $[0,T]$, and
therefore the pair $(A,B)$ is not controllable. This is a
contradiction and thus Lemma~\ref{claim1} is proved. \EOP

By standard homogeneity arguments Lemma~\ref{claim1}
\textcolor{black}{together with \rref{Vint}} imply
uniform exponential convergence of $V$ to zero along every
trajectory corresponding to a $(T,\mu)$-signal.
Theorem~\ref{thm-neutrally} is therefore proved. \EOP

\section{The double integrator}\label{sec:mainproof}

\subsection{Scalar control}

\textcolor{black}{This section addresses the same problem as
above for the double integrator. For this, let}
\begin{equation}\label{scandia}
A=\lp\ba{cc}0&1\\0&0\ea\rp,\ \ \ \textcolor{black}{B}=\lp\ba{c}0\\1\ea\rp,
\end{equation}
so that system \r{system} becomes
\be\label{DI}
\left\{
\ba{lcl}
\dot x_1&=&x_2,\\
\dot x_2&=&\alpha u,
\ea
\right.
\ee
while, by considering $K=(-k_1,-k_2)$, the closed-loop system \r{feedback} reads
\be\label{DIf}
\left\{
\ba{lcl}
\dot x_1&=&x_2,\\
\dot x_2&=&-\alpha (k_1 x_1+k_2 x_2)\,. \ea \right. \ee

\n Throughout the section we prove the following fact.
\begin{thm}\label{theo}
For every $0<\mu\leq T$ there exists a $(T,\mu)$-stabilizer for \r{DI}.
\end{thm}

\textcolor{black}{Again, we stress that, given any $T\geq \mu>0$,
the above result establishes the existence of a static linear feedback
that globally asymptotically stabilizes \rref{system} for the case that
$A$ and $B$ are given by \rref{scandia}, for \emph{all} $(T,\mu)$-signal $\alpha$.
Theorem \ref{theo} therefore gives a positive answer to Question $(Q2)$ for this particular case.}
\vspace{5mm}

\n \underline{Proof of Theorem~\ref{theo}.}  We first show how to
exploit the symmetries of system \r{DI} in order to identify a
one-parameter family of problems which are equivalent to the
research of a $(T,\mu)$-stabilizer.

\begin{lemma}\label{multi}
Let $\lb$ be a positive real number. Then \r{DI} admits a
$(T,\mu)$-stabilizer if and only if it admits a
$(T/\lb,\mu/\lb)$-stabilizer. More precisely, $(-k_1,-k_2)$ is a
$(T,\mu)$-stabilizer if and only if $(-\lb^2 k_1,-\lb k_2)$ is a
$(T/\lb, \mu/\lb)$-stabilizer.
\end{lemma}
\underline{Proof of Lemma~\ref{multi}.} Let $K=(-k_1,-k_2)$ be a
$(T,\mu)$-stabilizer and fix
 an arbitrary $\lb>0$. It suffices to prove that $(-\lb^2 k_1,-\lb k_2)$
is a $(T/\lb, \mu/\lb)$-stabilizer, the converse part of the statement being
equivalent (just making $1/\lb$ play the role of $\lb$).

Applying to $x(\cdot):=x(\cdot\,;t_0,x_0,K,\alpha)$ a time-rescaling and an anisotropic dilation, we define
\[
x_\lb(t)=\lp\ba{cc}1&0\\0&\lb\ea\rp x(\lb t)\,,\quad \forall
t\geq 0\,.
\]
Then
\brs
\frac{d}{dt}x_\lb(t)&=&\lb \lp\ba{cc}1&0\\0&\lb\ea\rp \dot x(\lb t)=
\lb \lp\ba{cc}1&0\\0&\lb\ea\rp \lp\ba{cc}0&1\\ -k_1\al(\lb
t)&-k_2\al(\lb t)\ea\rp x(\lb t)\\
&=&\lp\ba{cc}0&1\\ -\lb^2 k_1\al(\lb t)&-\lb k_2\al(\lb t)\ea\rp x_\lb(t)
=Ax_\lb(t)-\al(\lb t)b(\lb^2 k_1,\lb k_2) x_\lb(t)\,,
\ers
that is,
$$x_\lb(\cdot)=x(\cdot\,;t_0,\mathrm{Diag}(1,\lb)x_0,(-\lb^2 k_1,-\lb k_2),\alpha(\lb \cdot)).$$
It is clear that $\alpha(\lambda\cdot)$ is a $\lp T/\lb,
\mu/\lb\rp$-signal if and only if $\al(\cdot)$ is a $(T,\mu)$-signal.
Therefore, if $\beta$ is a class $\cKL$ function such that, for every
$\al\in\caG(T,\mu)$, $x_0\in\mR^2$ and $t_0\in\mRp$,
$$
\norm{x(t;t_0,x_0,K,\alpha)}\leq \beta(\norm{x_0},t-t_0)\,,\quad \forall t\geq t_0\,,
$$
then, for every  $\al_\lb\in\caG(T/\lb,\mu/\lb)$, $x_0\in\mR^2$ and $t\geq t_0\geq0$,
$$
\norm{x(t;t_0,x_0,(-\lb^2 k_1,-\lb k_2),\alpha_\lb(\cdot))}=
\norm{x(t;t_0,\mathrm{Diag}(1,1/\lb)x_0,K,\alpha(\cdot/\lb)}\leq
\beta(\max\{1,1/\lb\}\norm{x_0},t-t_0)\,.
$$
Since $(s,t)\mt \beta(\max\{1,1/\lb\}s,t)$ is a class $\cKL$ function,
the lemma is proved. \EOP

\vspace{5mm}
We prove below Theorem \ref{theo} by fixing a gain $K=(-k_1,-k_2)$
and showing that there exists $\lambda>0$  such that $(-\lambda^2
k_1,-\lambda k_2)$ is a $(T,\mu)$-stabilizer for \r{DI}. According
to Lemma~\ref{multi}, this is equivalent to proving that  there
exists $\lambda>0$  such that $K$ is a $(T/\lb,\mu/\lb)$-stabilizer.

The first, obliged, choice is on the sign of the entries of $K$: a
necessary (and sufficient) condition for the matrix $A+bK$ to be
Hurwitz is that $k_1,k_2>0$.

The choice of $K$ will be determined by the following request: we
ask that each matrix $A+\alpha b K$, with $\alpha\in[\mu/T,1]$ constant, has
real negative eigenvalues. Since the discriminant of
\[\det\lp\ba{cc}-\sigma&1\\ -\al k_1&-\sigma-\al k_2\ea\rp=
\sigma^2+\al k_2 \sigma+\al k_1\] is given by $\al(\al k_2^2-4
k_1)$, this sums up to imposing that
\be\label{k_int}
k_1<\frac{\mu }{4T}k_2^2.
\ee

We fix for the rest of the argument a positive $\rho<\mu/2T$
and take $K=(-\rho k^2/2,-k)$ for $k>0$ to be fixed later. The choice of $\rho$ and $K$ is such that
inequality \r{k_int} is automatically verified. For $\bar{\alpha}\in\{\mu/T,1\}$, define
$\xi^{\bar{\alpha}}_{\pm}$ as the roots of
$$
\xi^2+\bar{\alpha}k\xi+\bar{\alpha}\frac{\rho k^2}2=0,
$$
with $\xi^{\bar{\alpha}}_+<\xi^{\bar{\alpha}}_-$. In addition, let
$\xi^s_+$ and $\xi^s_-$ be given by
$$\xi^s_+:=-\frac k2(1+\sqrt{1-\rho}),\ \ \
\xi^s_-:=-\frac k2\lp1-\sqrt{1-\lp2-\frac\rho 2\rp\rho}\rp. $$ A simple
calculation shows that
\be\label{xi0}
\xi^s_+<\xi^1_+<\xi^{\mu/T}_+<\xi^{\mu/T}_-<\xi^1_-<\xi^s_-<0. 
\ee


We next define a set of half-lines of the plane. Let $D_1$ and $D_2$
be the half-lines defined, respectively, by $x_2=0,x_1<0$ and
$x_2=0,x_1>0$. Let, moreover,
$D^s_{\pm}$
and
$D^{\bar{\alpha}}_{\pm}$ be the half-lines of the open
upper half-plane defined, respectively, by the equations
\be\label{lines0}
D^s_{\pm}: \ \ x_2=\xi^s_{\pm} x_1, \
\ \ D^{\bar{\alpha}}_{\pm}: \ \ x_2=\xi^{\bar{\alpha}}_{\pm}x_1.
\ee

\begin{figure}[h!]
\begin{center}
\includegraphics[scale=0.8]{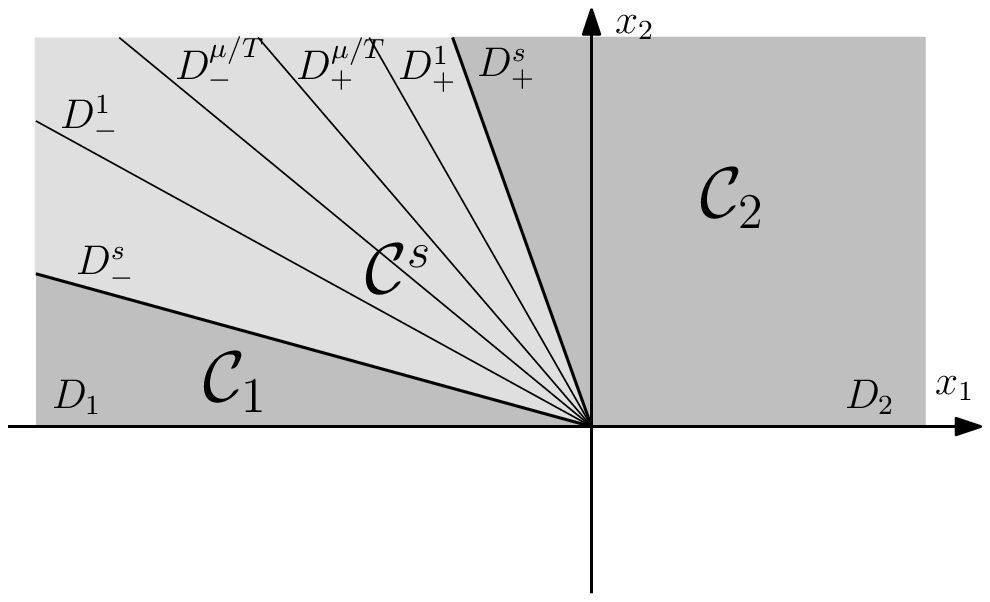}
\caption{The cones ${\cal{C}}_1$, ${\cal{C}}^s$, ${\cal{C}}_2$.}\label{halfplane}
\end{center}
\end{figure}

Finally, we define 
${\cal{C}}_1$, ${\cal{C}}^s$, and
${\cal{C}}_2$ as the closed cones contained in the upper half-plane and delimited, respectively, by $D_1$, $D^s_-$,
then $D^s_-$, $D^s_+$, and $D^s_+$, $D_2$: see Figure \ref{halfplane}.

A key step in the proof of Theorem~\ref{theo} is to evaluate
for how much time a trajectory of  \r{DIf} stays in each of the three cones.
\begin{lemma}\label{finite}
There exists a positive constant $C_1=C_1(\rho,T,\mu)$ ({\it i.e.}, only
depending on $\rho$, $T$, and $\mu$) such that, for every
$\lambda\geq
k$, every $(T/\lb,\mu/\lb)$-signal $\al$, every $x_*\in\R^2\setminus\{0\}$, and every $t_0\geq 0$, if an interval $I\subset \R$ is such that the
trajectory $t\mt x(t;t_0,x_*,K,\al)$ stays in ${\cal{C}}_1\cup {\cal{C}}_2$ for every $t\in I$ then the  length of $I$ is smaller than
$C_1/k$.
\end{lemma}

\underline{Proof of Lemma~\ref{finite}.}
Consider $\al\in L^\infty(\mRp,[0,1])$ and a trajectory $x(\cdot):=x(\cdot\,;t_0,x_*,K,\al)$ of \r{DIf} with $x_*\in\R^2\setminus\{0\}$ and $t_0\geq 0$. Let $I=[t_1,t_2]$
be a time-interval such that $x(t)$ belongs
to ${\cal{C}}_1\cup{\cal{C}}_2$ for every $t\in I$.

Using polar coordinates, $x(\cdot)$ can be represented
as $t\mt (r(t)\cos\th(t),r(t)\sin\th(t))$ and one has
$$
\dot\th=\frac{\dot x_1\sin\th-\dot x_2\cos\th}{r}=-\sin^2\th+\al\cos\th\lp\frac{\rho k^2}2\cos\th+k\sin\th\rp
$$
almost everywhere in $I$.
Therefore,
$$
\dot\th\leq -\al\lp \sin^2\th+k\cos\th\sin\th+\frac{\rho
k^2}2\cos^2\th\rp.$$
We next show that there exists a positive
constant $c(\rho)$ such that
\be\label{derPE}
\dot\th\leq-\alpha c(\rho)(\sin^2\th+k^2\cos^2\th),
\ee
almost everywhere in $I$. The claim can be proved by
taking $\tau:=\tan(\th)$ and noticing that there exists $c(\rho)>0$ small enough
such that if
$\tau<\xi^s_+$ or $\tau>\xi^s_-$ then
$$
\tau^2+k\tau+\frac{\rho k^2}2\geq c(\rho)(\tau^2+k^2).
$$

Define now
$$
F(\th)=\left\{\ba{ll}\arctan\lp \frac{\tan\th}{k}\rp&\mbox{if }0\leq \th<\pi/2,\\ \pi/2&\mbox{if }\th=\pi/2,\\ \arctan\lp \frac{\tan\th}{k}\rp+\pi&\mbox{if }\pi/2< \th\leq \pi.\ea\right.
$$
Notice that $F$ is a continuous reparameterization of the interval $[0,\pi]$.
From (\ref{derPE}), we get
\be\label{new}
\frac{d}{d t} (F(\th))\leq-c(\rho)k\alpha.
\ee

Therefore, $F(\th)$ is monotone non-increasing as long as $x(\cdot)$
stays in ${\cal{C}}_1\cup {\cal{C}}_2$.
If $\al$ is a
$(T/\lb,\mu/\lb)$-signal, then
\be\label{gap}
F\lp \theta\lp t+\frac T{\lb}\rp\rp -F(\th(t))\leq -\frac{c(\rho)\mu k}{\lb}
\ee
for every $t$ with $[t,t+T/\lb]\subset I$.
Let $l^{\max}$ be the largest 
integer such that
$I$ contains $l^{\max}$ disjoint sub-intervals of length $T/\lb$.
Then \r{gap} yields
\[ \pi \geq F(\th(t_1))-F(\th(t_2))\geq l^{\max}\frac{c(\rho)\mu k}{\lb} \]
and thus
\[l^{\max}\leq \frac{\lb\pi}{c(\rho)\mu k}.\]
Hence
\[t_2-t_1\leq (l^{\max}+1)\frac T\lb
\leq \frac{ T\pi}{c(\rho)\mu k}+\frac T\lb.\]
The lemma is concluded by recalling the hypothesis $\lb\geq k$.
\EOP

\vspace{5mm}
From now until the end of the proof, let $x_0\ne 0$ belong to  the cone
defined by $D_1$ and $D^s_+$, {\it i.e.}, to ${\cal{C}}_1\cup {\cal{C}}^s$.
Let, moreover,
$t_0\geq 0$, $\lambda\geq k$, $\al\in {\cal G}(T/\lb,\mu/\lb)$, and define $x(\cdot):=x(\cdot\,;t_0,x_0,K,\al)$.
Then  the following
alternative occurs for $x(\cdot)$:
\begin{description}
\item[$(i)$] for every $t\geq t_0$, $x(t)$ remains in ${\cal{C}}_1\cup {\cal{C}}^s$;
\item[$(ii)$] $x(\cdot)$ reaches $D^s_+$ in finite time.
\end{description}
In both cases, let $[t_0,t_1]$ be the time-interval  needed by
$x(\cdot)$ to reach ${\cal C}^s$. Recall that, by Lemma~\ref{finite}, one
has $0\leq t_1-t_0<C_1/k$.
We notice the following fact, which results from a trivial
computation.

\begin{lemma}\label{ff00}
The positive definite function $V(x):=x_1^2+2x_2^2/\rho k^2$,
evaluated along $x(\cdot)$, is non-increasing as long as $x(\cdot)$
remains in the fourth quadrant, {\it i.e.} $\{x\in\mathbb R^2 : x_1\leq 0, x_2\geq 0\}$.
\end{lemma}

\n The following lemma provides an exponential decay result for trajectories staying in ${\cal C}^s$ (in particular, for trajectories satisfying ($i$)).
\begin{lemma}\label{ff01}
There exist
$C_2=C_2(\textcolor{black}{\rho,k},T,\mu)>0$ and $\gamma=\gamma(\rho,T,\mu)>0$
such that,
for every $t\geq t_1$ such that $x(\cdot)$ stays in ${\cal C}^s$ along the interval $[t_1,t]$, it holds that
 \be\label{ff11} |x(t)|\leq
C_2 e^{-k\gamma (t-t_1)}|x_0|.
\ee
Moreover, $C_2=O(k)$ as a function of $k$ as $k$ tends to infinity.
\end{lemma}
\underline{Proof of Lemma~\ref{ff01}.}
Let $t$ be as in the statement of the lemma.
We deduce that, for $\tau\in (t_1,t)$,
$$
\dot x_2(\tau)=-k\al(\tau)w(\tau)x_2(\tau),
$$
where $w(\cdot)$ is a continuous function verifying
$$
0<1+\frac{\rho k}{2\xi^s_-}\leq w(\tau)\leq 1+\frac{\rho k}{2\xi^s_+}.
$$
Notice that the bounds on $w$ do not depend on $k$, due to the definition of $\xi^s_\pm$. We are back to the (PE) one-dimensional case studied in the
introduction. We deduce that there exist $C_*=C_*(\rho,T,\mu)>0$ and
$\gamma=\gamma(\rho,T,\mu)>0$ such that
\be\label{inter00}
x_2(t)\leq C_* e^{-k\gamma(t-t_1)}x_2(t_1).
\ee

Notice now that, for every $x\in {\cal C}^s$,
$$|x|\leq x_2\sqrt{1+\frac1{(\xi_+^s)^2}},\ \ \ \ \ \ x_2\leq k\sqrt{\frac\rho2}\sqrt{V(x)},\ \ \ \ \ \ \sqrt{V(x)}\leq \max\left\{1,\textcolor{black}{\frac 1 k\sqrt{\frac 2\rho}}\right\}|x|.$$

The proof of the lemma is concluded by
recalling  that $V(x(t_1))\leq V(x_0)$ and
plugging the above estimates in \r{inter00}.
\EOP

\vspace{5mm}
Let us now
establish a lower bound on the time needed to go across
 the cone
${\cal C}^s$.
\begin{lemma}\label{final0}
There exists $\lambda_0>0$ such that
if $\lb\geq\lb_0$ and $x_0\in D^s_-$ then $x(\cdot)$ stays in ${\cal C}^s$ for all $t\in [t_0,t_0+1]$.
\end{lemma}
\underline{Proof of Lemma~\ref{final0}.}
Fix $x_0\in D^s_-$.
Reasoning by contradiction (and exploiting the homogeneity of \r{DIf}) we assume that there exist a strictly
increasing unbounded sequence $\{\lb_i\}_{i\in \mathbb N}$
and a sequence $\{\alpha_i\}_{i\in \mathbb N}\subset L^\infty(\mRp,[0,1])$ with $\alpha_i\in
\caG(T/\lb_i,\mu/\lb_i)$ such that each
$x_i(\cdot):=x(\cdot\,;t_0,x_{0},K,\al_i)$, $i\in \mathbb N$,
reaches $D^s_+$ in time smaller than one. Because of the sequential
weak-$\star$  compactness of $L^\infty(\mRp,[0,1])$, there
exists $\alpha_\star\in L^\infty(\mRp,[0,1])$ and a
subsequence of $\{\lb_i\}_{i\in \mathbb N}$ (still denoted by
$\{\lb_i\}_{i\in \mathbb N}$) such that $\{\alpha_i\}_{i\in \mathbb N}$ converges
weakly-$\star$ to $\alpha_\star$.

By Proposition~\ref{technic},
we deduce that the sequence $\{x_i\}_{i\in \mathbb N}$ converges, uniformly
on compact time-intervals, and in particular on $[t_0,t_0+1]$,
to $x_{\star}(\cdot):=x(\cdot\,;t_0,x_0,K,\al_\star)$. Since, for every
$i\in \mathbb N$, $x_i(\cdot)$ reaches $D^s_+$ in time smaller than
one, we deduce that $x_{\star}(t_0+t_{\star})\in D^s_+$ for some
$t_{\star}\in (0,1]$.

Let us show that $\al_{\star}\geq \mu/T$ almost everywhere on
$\mRp$. For every interval $J\subset\mRp$ of finite length $\ell>0$,
apply  \r{weakstar} to the characteristic function
of $J$.
Since each
$\alpha_{\lb_i}$ is a $(T/{\lb_{i}},\mu/{\lb_{i}})$-signal, it
follows that
\[
\frac 1 \ell \int_J \alpha_\star(s)ds=\lim_{i\to\infty}
\frac 1 \ell \int_J \alpha_{i}(s)ds \geq \liminf_{i\to\infty} \frac{\mu}{\ell \lb_{i}}
{\cal I}\left(\frac{\ell\lb_{i}}{T}\right)=\frac{\mu}{T}\,,
\]
where ${\cal I}(\cdot)$ denotes the integer part.
Recall that since $\al_\star$ is measurable and bounded (actually, $L^1$ would be enough),
almost every $t>0$ is a Lebesgue point for $\al_\star$, {\em i.e.}, the limit
\[\lim_{\eps\to0+}\frac 1 {2\eps} \int_{t-\eps}^{t+\eps} \alpha_\star(s)ds\]
exists and is equal to $\al_\star(t)$ (see, for instance,
\cite{leb_points}). We conclude that, as claimed,
$\alpha_\star(t)\geq \mu/T$ almost everywhere.

Therefore,
$x_{\star}(\cdot)$ is actually a trajectory of
the switching system
\be\label{swi0}
\dot x =\left[u(A+bK^\top)+(1-u)\lp A+\frac{\mu}TbK^\top\rp\right]x,
\ee
where $u$ is a measurable
function defined on $\mRp$ and taking values in $[0,1]$ (see
\cite{BOSCAIN2002} and references therein for more on switching
systems). According to the taxonomy and the results in \cite[page 93]{BOSCAIN2002}, \r{swi0}
is a switching system of type
$(RR.2.2.A)$ and, as a consequence, the curve
$x_{\star}$ stays below the trajectory $\bar{x}$ of \r{swi0} starting from
$x_0$ at time $t_0$ and corresponding to the input $u\equiv 0$. 
Since $\bar x$ converges to zero in the cone delimited by $D^s_-$ and $D^{\mu/T}_-$, $x_{\star}$ must stay in the same
cone, contradicting the fact that $x_{\star}$ reaches
$D^s_+$ in finite time. Lemma~\ref{final0} is proved. \EOP

\vspace{5mm}
Let us now focus on the behavior of trajectories exiting ${\cal C}^s$ or, equivalently,
such that $x_0\in D^s_+$.
\begin{lemma}\label{c2}
There exists $\lambda_1\geq k$ such that
if $\lb\geq \lb_1$ and $x_0\in D^s_+$ then \textcolor{black}{there exists a finite time
$t_f>0$ such that} $x(\cdot)$ satisfies
\be\label{c22}
|x(t_f)|\leq|x_0|,
\ee
with $x(t_f)\in D_2$ and $x(t)\in{\cal C}_2$ for all $t\in [t_0,t_f]$.
\end{lemma}

\underline{Proof of Lemma~\ref{c2}.} Fix $x_0\in D^s_+$. We reason again by
contradiction and follow the same procedure as in Lemma~\ref{final0}.
This is possible since, according to
Lemma~\ref{ff01}, the time needed by any trajectory of \r{DIf} to
go across in ${\cal {C}}_2$ is bounded (uniformly with respect to
$\lb\geq k$ and $\alpha\in\caG(T/\lb,\mu/\lb)$) by $C_1/k$.
We obtain that,
for some $t_0<t'_f\leq C_1/k$, the limit trajectory $x_\star$
is contained in ${\cal C}_2$ on $[t_0,t_f']$, reaches $D_2$ at time $t_f'$, and satisfies  $|x_{\star}(t'_f)|\geq|x_0|$.
 According
to \cite{BOSCAIN2002}, the trajectory $x_{\star}$ is, inside ${\cal
{C}}_2$, below the integral curve $x_{\mu/T}$ of $\dot x=
(A+\frac{\mu}TbK^\top)x$ with initial
condition $x_{\mu/T}(t_0)=x_0$. That means in particular that
$|x_{\star}(t'_f)|\leq |x_{\mu/T}(t''_f)|$, where $t''_f$ is the first time larger than $t_0$ such that $x_{\mu/T}(t''_f)\in D_2$. However, a lengthy but straightforward computation
shows that $|x_{\mu/T}(t''_f)|<|x_0|$ and we reach a
contradiction. \EOP

\vspace{5mm}
Define $\lb_*=\max\{\lb_0,\lb_1\}$ where $\lb_0$ and $\lb_1$ are the quantities appearing in the statements of Lemmas~\ref{final0} and \ref{c2}.
Whenever $x(\cdot)$ satisfies ($ii$) we let $(t_1\leq)t_2<t_f$ be such that
$x(t)\in{\cal C}^s$ for $t\in[t_1,t_2]$, $x(t)\in{\cal C}_2$ for $t\in[t_2,t_f]$, and $x(t_f)\in D_2$.

As a final technical result for the completion of the argument,
we need the following lemma.
\begin{lemma}\label{ouf0}
There exist $k_*=k_*(\rho,T,\mu)>0$ and
$\gamma_*=\gamma_*(\rho,T,\mu)>0$
such that
if $x(\cdot)$ satisfies $(ii)$,
$t_f-t_0\geq 1$, $k\geq k_*$, and $\lb\geq \lb_*$,
then
\be\label{ouf}
|x(t_f)|\leq
\frac12 e^{-k\gamma_* (t_f-t_0)}|x_0|.
\ee
\end{lemma}

\underline{Proof of Lemma~\ref{ouf0}.}
Applying estimates (\ref{ff11}) and (\ref{c22}) we get that
\[
|x(t_f)|\leq|x(t_2)|\leq C_2 e^{-k\gamma (t_2-t_1)}|x_0|\,.
\]
Moreover, according to Lemma~\ref{finite} and the hypothesis $t_f-t_0\geq 1$,
\[
|x(t_f)|\leq C_2 e^{2\gamma C_1} e^{-k\gamma (t_f-t_0)}|x_0|\leq C_2 e^{2\gamma C_1-k\gamma/2} e^{-k\gamma (t_f-t_0)/2}|x_0|\,.
\]
The lemma is proved by taking $\gamma_*=\gamma/2$ and
$k_*$ large enough in order to have  $C_2 e^{2\gamma C_1- k\gamma/2}\leq 1/2$. Such a $k_*$ does exist because, in view of Lemmas \ref{finite} and \ref{ff01}, neither $C_1$ nor $\gamma$  depend on $k$, while $C_2=O(k)$ as $k$ tends to infinity.
 \EOP

\vspace{5mm}
We have developed enough tools to conclude the proof of Theorem~\ref{theo}.
Take $k\geq k_*$ and $\lb\geq \lb_*$.
Let $\xi_0\in\R^2\setminus\{0\}$, $t_0\geq 0$, $\al\in{\cal G}(T/\lb,\mu/\lb)$, and denote by $\Tau$ the set of
times such that $\xi(\cdot)=x(\cdot\,;t_0,\xi_0,\al)$ belongs to the $x_1$-axis.
Choose one representative for every connected component of $\Tau$.
Denote by $\Tau'$ the set of such representatives and by $j\in\N\cup\{+\infty\}$ its cardinality. By monotonically enumerating the elements of $\Tau'$, we have
$\Tau'=\{\tau_i\mid 0\leq i<j\}$ with $\tau_{i-1}< \tau_{i}$ for $0<i<j$.
Let, moreover, $\tau_{-1}=t_0$ and $\tau_{j}=+\infty$.
Since,  for every $-1\leq i<j$, either $\xi(\cdot)$ or $-\xi(\cdot)$, restricted to $(\tau_{i},\tau_{i+1})$, is contained in the upper half-plane, then the previously established estimates apply to it.

Take
$-1\leq i<j$ and $t\in[\tau_i,\tau_{i+1})$.
Then there exists $C_3=C_3(K,T,\mu)>1$ such that
\be\label{kind}
|\xi(t)|\leq C_3 e^{-k\gamma_* (t-\tau_{i})}|\xi(\tau_{i})|.
\ee
Indeed, let
\brs
s_1&=&\sup\{s\in [\tau_i,t]\mid \xi([\tau_i,s])\subset {\cal C}_1\mbox{ or }-\xi([\tau_i,s])\subset {\cal C}_1\},\\
s_2&=&\inf\{s\in [\tau_i,t]\mid \xi([s,t])\subset {\cal C}_2\mbox{ or }-\xi([s,t])\subset {\cal C}_2\},
\ers
and notice that $\max\{s_1-\tau_i,t-s_2\}< C_1/k$. Then \r{kind} immediately follows from
Lemma~\ref{ff01}.

In particular, if $j>0$
then
$$ |\xi(\tau_{0})|\leq
C_3 e^{-k\gamma_* (\tau_{0}-t_{0})}|\xi_0|.
$$
Moreover, for every $0<i<j$, Lemma~\ref{ouf0} yields
$$ |\xi(\tau_{i})|\leq
\frac12 e^{-k\gamma_* (\tau_{i}-\tau_{i-1})}|\xi(\tau_{i-1})|,
$$
which implies, by recurrence, that
$$ |\xi(\tau_{i})|\leq
\frac12 e^{-k\gamma_* (\tau_{i}-\tau_{0})}|\xi(\tau_{0})|.
$$
Therefore, independently of $j$, for every $-1\leq i<j$ we have
$$ |\xi(\tau_{i})|\leq C_3 e^{-k\gamma_* (\tau_{i}-t_{0})}|\xi_{0}|. $$
Applying again \r{kind} we obtain that for every $t\in[t_0,+\infty)=\bigcup_{i=1}^{j-1}[\tau_{i-1},\tau_i)$,
$$ |\xi(t)|\leq C_3^2 e^{-k\gamma_* (t-t_{0})}|\xi_{0}|\,, $$
which proves Theorem~\ref{theo}.
\EOP

\vspace{5mm}
Let us go back to the natural question posed in $(Q1)$, that is, whether the choice of a stabilizer
$K$ can be made independently of $T$ and $\mu$.
Unlike the case when $A$
is neutrally stable, we prove below that the answer is negative when
the double integrator is considered.

\begin{prop}\label{prop1}
For every $K\in\R^2$, there exist $T\geq\mu>0$ such that $K$ is not
a $(T,\mu)$-stabilizer for system \r{DI}.
\end{prop}

\begin{figure}[h!]
\begin{center}
\includegraphics[scale=0.4]{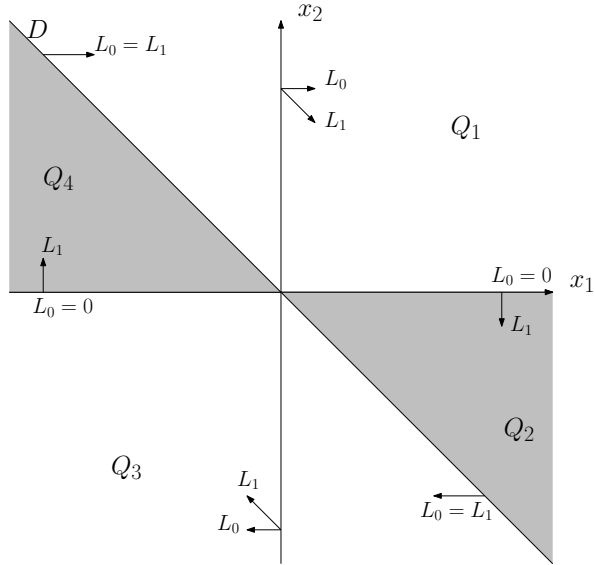}
\end{center}
\caption{Partition of the state-space into the four regions
$Q_1,Q_2,Q_3,Q_4$.}\label{quadrants}
\end{figure}

\underline{Proof of Proposition~\ref{prop1}.} Let
$K=(-k_1,-k_2)\in\R^2$. If $k_1\leq 0$ or $k_2\leq0$ then $A+b K$ is
not Hurwitz and thus $K$ is not  a $(T,\mu)$-stabilizer for system
\r{DI}, whatever $\mu$ and $T$. Let now $k_1,k_2>0$. Among the possible configurations of the right-hand side of
\r{DIf}, we focus our attention on the linear vector fields $L_0$ and
$L_1$, corresponding to $\alpha(t)=0$ and to $\alpha(t)=1$
respectively, that is,
$$
L_0(x):=A\,x\,,\quad L_1(x):=\lp \ba{cc}0&1\\ -k_1& -k_2\ea \rp x\,.
$$
Consider the set where $L_0$ and $L_1$ are collinear, {\em i.e.},
the union of the axis $\{x_2=0\}$ and the line $D$ defined by
$$
D:=\{(x_1,x_2)\in\mR^2\mid x_2=-(k_1/k_2)x_1\}\,.
$$
We denote by $Q_1,Q_2,Q_3,$ and $Q_4$ the four regions of the plane
delimited by these two lines and defined as follows
\begin{eqnarray*}
Q_1&:=& \{(x_1,x_2)\in\mR^2\mid x_2\geq -(k_1/k_2)x_1\,, x_2> 0\},\\
Q_2&:=& \{(x_1,x_2)\in\mR^2\mid x_2> -(k_1/k_2)x_1\,, x_2\leq 0\},\\
Q_3&:=& \{(x_1,x_2)\in\mR^2\mid x_2\leq -(k_1/k_2)x_1\,, x_2< 0\},\\
Q_4&:=& \{(x_1,x_2)\in\mR^2\mid x_2< -(k_1/k_2)x_1\,, x_2\geq 0\}.
\end{eqnarray*}
Roughly speaking, on $Q_2\cup Q_4$ the vector field $L_1$ points
``more outwards'' than $L_0$  (with respect to the origin). Therefore,
on $Q_2\cup Q_4$ the excitation ($\al=1$) is not helping the stability
towards the origin (see Figure~\ref{quadrants}).

We look for a ``destabilizing" $\al$ in the form $\zeta(x(\cdot))$, where
$x(\cdot)$ is an unbounded solution corresponding to the feedback
\[
\zeta(x)=\left\{\ba{ll}
1&\mbox{if }x\in Q_2\cup Q_4,\\
\mu/T&\mbox{if }x\in Q_1\cup Q_3,
\ea\right.
\]
$\mu$ and $T$ to be chosen. Clearly, any $\al$ obtained in such
a way is a $(T,\mu)$-signal. (It should be stressed that
although $\zeta$ is a discontinuous feedback, its solutions are anyway well-defined
because all trajectories of \r{DI} are rotating around the origin).

Let $\bar x$ be the first intersection with $D$ of the trajectory of
\r{DIf} starting from $(-1,0)$ and corresponding to $\al\equiv 1$.
Notice that there exists
$\nu>0$ small enough
such that
 the trajectory of \r{DIf} starting from
$\bar x$ and corresponding to $\al\equiv \nu$ crosses the $x_1$-axis
at a point $(\xi,0)$ such that $\xi>1$.  This is because
as
$\nu\to 0$ such trajectories converge uniformly on compact intervals to the horizontal curve $t\mt e^{t A}\bar x$.

Fix $T\geq \mu>0$ such that
$\mu/T\geq \nu$. Then the trajectory $x(\cdot)$ of
\r{DIf} starting from $(-1,0)$ and corresponding to the feedback $\zeta$ first leaves the upper half-plane through $[\xi,+\infty)\times\{0\}$.
Then, by symmetry and
homogeneity of the system, $|x(t)|$ goes to infinity as $t\to\infty$.
Thus, $\al(t)=\zeta(x(t))$ is the required destabilizing signal.
\EOP

\subsection{Non-scalar control}
It makes sense to consider the stability properties of the system whose linear dynamics is the same as that of the double integrator, but which has a
different controlled part. That is, we study system
 \be\label{multiDI}
\dot x=\lp\ba{cc}0&1\\0&0\ea\rp x+\al B u,\ \ \ \ \ x\in\R^2,\ \ \ \ \ u\in\R^m,
\ee
where $B$ is a general $2\times m$ matrix such that the pair $(A,B)$ is controllable
($A$ denotes the nilpotent $2\times2$ matrix appearing in \r{multiDI}).
Since the image of $A$ is one-dimensional,
the controllability of $(A,B)$ implies the existence of  a column $b$ of $B$ such that $(A,b)$ is controllable.
Moreover, by a linear change of coordinates,
we can transform $(A,b)$ into its Brunovsky normal form, that is, into the
matrix $A$ and the vector $b=(0 , 1)^\top$ as in \r{scandia}.
Therefore Theorem~\ref{theo} guarantees that for every pair of positive constants $\mu\leq T$ there exists a $(T,\mu)$-stabilizer for \r{system}.

The difference with respect to the scalar case is that
Proposition~\ref{prop1}
is not valid anymore when the rank of $B$ equals two, that is, the answer to question $(Q1)$ becomes positive.

\begin{prop}\label{q1yes}
If $\mathrm{rank}(B)=2$, then
there exists $K$ of size $m\times 2$ such that
for every  $T\geq\mu>0$
the gain $K$ is a
$(T,\mu)$-stabilizer for system \r{multiDI}.
\end{prop}

\underline{Proof of Proposition~\ref{prop1}.} Let us remark that, up to a reparameterization of $\R^m$ of the type $u'=M u$, $M$ invertible, we can rewrite $B$ as $(\Id\ {\bf 0})$, where $\Id$ is the $2\times 2$ identity matrix and all the entries of the $2\times(m-2)$ matrix ${\bf 0}$ are null.
That is, \r{multiDI} is equivalent to
 \be\label{multiDI2}
\dot x=A x+\al u,\ \ \ \ \ x\in\R^2,\ \ \ \ \ u\in\R^2.
\ee
Fix $k>0$ and take $u=-k x$. If $t\mt x(t)$ is a solution of
\[\dot x=(A -\al k\Id)x,\]
then $t\mt y(t):=e^{-At}x(t)$ satisfies $\dot y=-k\al y$. Therefore
$$
y(t)= \exp\left(-k\int_0^t\al(\tau)d\tau\right) y(0)\,,
$$
which implies that
$$
\norm{x(t)}\leq \norm{e^{A t}}\exp\left(-k\int_0^t\al(\tau)d\tau\right)\norm{x(0)}.
$$
Since $e^{A t}=\Id+tA$ is linear in $t$,
we have proved that for every $T\geq\mu>0$ and for every $k>0$ the gain $K=-k\Id$ is a
$(T,\mu)$-stabilizer for system \r{multiDI2}.
\EOP

%

\section{Appendix}

We first provide a simple result used several times in the
paper.

\begin{prop}\label{technic}
Consider system \rref{feedback}, where $A,B,K$ are matrices of size
$n\times n$, $n\times m$ and $m\times n$ respectively, and $\alpha$
is a $(T, \mu)$-signal. Take a sequence $x_{0i}$ of norm-one
vectors converging to $x_{0*}$ and a bounded sequence $\alpha_i$ in
$L^\infty({\mathbb R})$ which converges weakly-$\star$ (in
$L^\infty({\mathbb R})$) to a measurable function $\alpha_\star$.
Then the sequence
$\{x_{i}(\cdot):=x(\cdot,0,x_{0i},K,\alpha_i)\}_{i\in\mathbb N}$ converges,
uniformly on compact time intervals, to
$x_{\star}(\cdot):=x(\cdot,0,x_{0\star},K,\alpha_\star)$ as $i$
tends to infinity.
\end{prop}

\underline{Proof of Proposition~\ref{technic}.} Recall that the
weakly-$\star$ convergence of $\al_i$ to $\al_\star$ means that, for
every $\varphi\in L^1(\mRp,\mR)$, it holds that
\[
\lim_{i\to\infty} \int_0^\infty \alpha_{i}(s)\varphi(s)ds =
\int_0^\infty \alpha_\star(s)\varphi(s)ds\,.
\]
Taking as $\varphi$ the characteristic function of an arbitrary interval
of length $T$
shows that $\alpha_\star$ is a
$(T,\mu)$-signal.
Moreover, the norm of $x_{0*}$ is equal to one. For $i\in \mathbb N$, set
$$
e_i(\cdot)=x_{i}(\cdot)-x_{\star}(\cdot), \ \
e_{0i}=x_{0i}-x_{0*},$$
$$
A_i(\cdot)=A+\alpha_i(\cdot)BK,
$$
and let $\phi_i(\cdot,\cdot)$ be the fundamental solution of $\dot
x=A_i(t)x$. Integrating the differential equation verified by $e_i$,
one gets, for $t\geq 0$,
\[
e_i(t)=\int_0^tA_i(s)e_i(s)ds+h_i(t),
\]
where
$h_i(t)=\phi_i(t,0)e_{0i}+\int_0^t(\alpha_i(s)-\alpha_\star(s))BKx_{\star}(s)ds$.
Note that the functions $h_i$ are uniformly bounded over compact
time-intervals and the sequence they define converges point-wise to
zero as $i$ tends to infinity. Therefore, by combining Gronwall
Lemma and the bounded convergence theorem, it follows that the
sequence $\{e_i(\cdot)\}_{i\in\mathbb N}$ converges point-wise to
zero as $i$ tends to infinity. The uniform convergence on compact
time intervals results from the above combined with Ascoli theorem.
\EOP


\label{sec:concl}

\newcommand{\auth}[1]{{\rm #1}}
\newcommand{\tit}[1]{{``#1''}}
\newcommand{\titl}[1]{{\em #1}}
\newcommand{\jou}[1]{{\em #1}}
\newcommand{\vol}[1]{vol.~{#1}}
\newcommand{\pp}[1]{pp.~#1}

\newcommand{\SortNoop}[1]{}
  \def\nesic{Nesi\'{c}\,} %
  \def\astrom{{\SortNoop{As}\AA}str{\"{o}}m\,}\let\c=\cedille

\end{document}